\magnification=\magstep1
\def\R{I \!\! R}
\def\N{I \!\! N}
\def\BbbR{\R}
\def\BbbN{\N}
\def\bull{\vrule height .9ex width .8ex depth -.1ex}
{\footline={$\triangle$ \hfill}
\null
\vskip 1.15truein
\centerline{Some New Asymptotic Properties for the Zeros of}
\centerline{Jacobi, Laguerre and Hermite Polynomials}
\smallskip
\centerline{by}
\smallskip
\hbox to \hsize{
\hss
\llap{\vtop{\hbox{Holger Dette}\hbox{Universit\"at G\"ottingen}}}
	  \quad and\quad
\rlap{\vtop{\hbox{William J. Studden}\hbox{Purdue University}}}
\hss
}
\medskip
\centerline{Technical Report \# 93-5}
\vskip 2in
\centerline{Department of Statistics}
\centerline{Purdue University}
\bigskip
\centerline{January, 1993}
\centerline{July, 1993}
\vfill\break
}
\pageno=1

\centerline{SOME NEW ASYMPTOTIC PROPERTIES FOR THE ZEROS OF}
\centerline{JACOBI, LAGUERRE AND HERMITE POLYNOMIALS}

   \footnote{}{\hskip-.28truein$^*$Research supported in part by the 
	 Deutsche Forschungsgemeinschaft.}
   \footnote{}{\hskip-.28truein$^{**}$Research supported in part by NSF
	 Grant DMS 9101730.}
   \footnote{}{\hskip-.3truein AMS Subject Classification: 33C45.}
   \footnote{}{\hskip-.3truein Keywords and phrases: Jacobi
	 polynomials, Laguerre polynomials, Hermite polynomials,}
   \footnote{}{\hskip1.5truein $\!\!$asymptotic zero distribution,
	 chain sequences}
\smallskip
\centerline{by}
\smallskip
\hbox to \hsize{
\hss
\llap{\vtop{\hbox{Holger Dette$^*$}\hbox{Institut f\"ur Mathematische
       Stochastik}\hbox{Technische Universit\"at Dresden}
\hbox{Mommsenstr.~13}
	  \hbox{01062 Dresden}\hbox{GERMANY}}}
	  \quad and\quad
\rlap{\vtop{\hbox{William J. Studden$^{**}$}\hbox{Department of Statistics}
	  \hbox{Purdue University}\hbox{1399 Mathematical Sciences Bldg.}
	  \hbox{West Lafayette, IN\ \ 47907-1399}\hbox{U.S.A.}}}
\hss
}
\bigskip
\bigskip
\bigskip
\centerline{ABSTRACT}

\bigskip
For the generalized Jacobi, Laguerre and Hermite polynomials
$P_n^{(\alpha_n, \beta_n)} (x), L_n^{(\alpha_n)} (x),$\break
$H_n^{(\gamma_n)} (x)$ the limit distributions of the zeros are found,
when the sequences $\alpha_n$ or $\beta_n$ tend to infinity with a
larger order than $n$.
The derivation uses special properties of the sequences in the
corresponding recurrence formulae.
The results are used to give second order approximations for the
largest and smallest zero which improve (and generalize) the limit
statements in a paper of Moak, Saff and Varga [11].

\bigskip
\bigskip
\bigskip
\medskip

\baselineskip=15.8pt

\noindent {\bf 1.\ \ Introduction.}
The aim of this paper is to give some new asymptotic results for the
zeros of the classical orthogonal polynomials.
To be precise let $P_n^{(\alpha, \beta)} (x), L_n^{(\alpha)} (x)$ and
$H_n^{(\gamma)} (x)$ denote the Jacobi, Laguerre and Hermite polynomials
orthogonal with respect to the measures $(1 - x)^\alpha (1+x)^\beta dx
\ (\alpha, \beta > -1, x \in [-1, 1]), x^\alpha e^{-x} dx \ (\alpha > -1,
x \in [0, \infty))$ and $|x|^{2\gamma} e^{-x^2} dx \ (\gamma > -{1\over
2}, x\in (-\infty, \infty))$, respectively (see Szeg\"o [15]).
In the following we will allow the weight functions to depend on the
degree of the polynomials and consider the ``generalized'' classical
orthogonal polynomials $P_n^{(\alpha_n, \beta_n)} (x),
L_n^{(\alpha_n)} (x)$ and $H_n^{(\gamma_n)} (x)$ where $(\alpha_n)_{n
\in {\BbbN}}, (\beta_n)_{n\in {\BbbN}}$ and $(\gamma_n)_{n\in {\BbbN}}$
are sequences of real numbers $(\alpha_n, \beta_n > -1,
\gamma_n > -{1\over 2})$.
These polynomials are, in general, not orthogonal on the
corresponding intervals, which reduces the number of methods
for investigating  their properties.
If $\alpha_n, \beta_n$ are linear functions of $n$ various asymptotic
results for the generalized Jacobi polynomials can be found in
papers of Moak, Saff and Varga  [11], Gonchar and Rakhmanov [8],
Mhaskar and Saff [10], Gawronski and Shawyer [7], Chen and Ismail
[2] and Ismail and Li [9].
A different approach (fixing the degree of the polynomials and
varying the parameters) was discussed in Elbert and Laforgia
[4].
Generalized Laguerre and Hermite polynomials with parameters depending
linearly on $n$ have been studied in Chen and Ismail [2] and
Gawronski [6].
In this paper we will investigate the asymptotic behaviour of the
zeros of the generalized classical orthogonal polynomials when the
sequences of the parameters tend to infinity with a larger order than
$n$.

In Section 2 we find the limit distribution for the (suitable
standardized) zeros of the generalized Jacobi polynomials while
Section 3 states similar results for the generalized Laguerre and
Hermite polynomials. It turns out that there exist essentially three
types of limit distributions (Corollary 2.3, 2.4, 2.5  and 2.6) depending
on the order of $\alpha_n$, $\beta_n$ and $\gamma_n$.
Section 4 deals with second order approximations for the largest and
smallest zero of the classical orthogonal polynomials which improve
the statements of  [11] for the Jacobi case.
The results of this paper are based on certain characterizing
properties of the chain sequences of the Jacobi, Laguerre and
Hermite measure which have recently been established in
a paper of Dette and Studden [3]. This approach allows a
 very simple derivation of asymptotic properties for the
 zeros of the corresponding  ``generalized'' orthogonal polynomials.

\bigskip
\bigskip
\bigskip\noindent
{\bf 2.\ \ The asymptotic distribution of the zeros of the
generalized Jacobi polynomials.}
Let $\alpha_n, \beta_n > -1$ for all $n\in {\BbbN}$, then the zeros
of $P_n^{(\alpha_n, \beta_n)} (x)$ are all simple and located in the
interval $(-1, 1)$.
For $\xi \in {\BbbR}$ we define
$$
N^{(\alpha_n, \beta_n)} (\xi)~\colon = ~\# \{x~|~ P_n^{(\alpha_n,
   \beta_n)} (x) = 0,~ x\le \xi\}
\leqno(2.1)
$$
as the number of zeros of $P_n^{(\alpha_n, \beta_n)} (x)$ that are
less or equal than $\xi$ and $\mu^{(\alpha_n, \beta_n)}$ as the
(discrete) uniform distribution on the set $\{x|P_n^{(\alpha_n,
\beta_n)} (x) = 0\}$, i.e. $\mu^{(\alpha_n, \beta_n)} (\xi) = {1\over
n} N^{(\alpha_n, \beta_n)} (\xi)$.
For every probability measure on
the interval $[-1, 1]$ let $\Phi_\mu (z)$ denote
the Stieltjes transform of $\mu$ with corresponding
continued fraction expansion
\newbox\DigitBox \setbox\DigitBox=\hbox{0}
\newdimen\DigitWidth \DigitWidth=\wd\DigitBox
\catcode`@=\active \def@{\kern\DigitWidth}
$$
\leqalignno{
\Phi_\mu (z)~\colon &= ~\int_{-1}^1 {d\mu (x)\over z - x} ~=~
   {\lower6pt\hbox{@@@1@@@$|$}\over \raise6pt\hbox{$|z@+@1$}} -
   {\lower6pt\hbox{$2\zeta_1|$}\over \raise6pt\hbox{$|$@1@}} -
   {\lower6pt\hbox{$@2\zeta_2@|$}\over \raise6pt\hbox{$|z+1$}} -
   \ldots & (2.2)\cr
   & \cr
&= ~{\lower6pt\hbox{@@@@@1@@@@@@$|$}\over \raise6pt\hbox{$|z@+@1 -
   2\zeta_1$}} - {\lower6pt\hbox{@@@@@@$4\zeta_1 \zeta_2@@@@@|$}\over
   \raise6pt\hbox{$|z+1 - 2 (\zeta_2 + \zeta_3)$}} -
   {\lower6pt\hbox{@@@@@@$4 \zeta_3 \zeta_4@@@@@|$}\over
   \raise6pt\hbox{$|z+1 - 2 (\zeta_4 + \zeta_5)$}} -\ldots\cr
}
$$
where $\zeta_1 = p_1, \zeta_j = (1 - p_{j-1}) p_j$ $(j\ge 2)$ and the
quantities $p_j$ in the chain sequences $(\zeta_j)_{j\in {\BbbN}}$
are bijective continuous functions of the moments of the measure $\mu$
(see Perron [13], Wall [16] and Skibinsky [14]).
The following result characterizes the uniform distribution
$\mu^{(\alpha_n, \beta_n)}$ in terms of the quantities
$p_j$ and is an immediate
consequence from Lemma 2.1, Lemma 2.2 and Theorem 3.1 of Dette and
Studden [3].

\bigskip\noindent
{\bf Proposition 2.1.}~~
The uniform distribution $\mu^{(\alpha_n, \beta_n)}$ on the set $\{x|
P_n^{(\alpha_n, \beta_n)} (x) = 0\}$ is uniquely determined by the
chain sequence
$$
\leqalignno{
p_{2i}^{(n)} &~=~ {n-i\over 2 (n-i) + 1 + \alpha_n + \beta_n}\cr
&\hskip 2.5truein i = 1, \ldots, n &(2.3)\cr
p_{2i-1}^{(n)} &~=~ {\beta_n + n-i + 1\over 2 (n - i + 1) + \alpha_n +
   \beta_n}\cr
}
$$
in the
corresponding continued fraction expansion (2.2).

\bigskip\noindent
{\bf Theorem 2.2.}~~
Let $N^{(\alpha_n, \beta_n)} (\xi)$ be defined in (2.1) and assume
that there exists sequences $(\delta_n)_{n\in {\BbbN}}$ and
$(\varepsilon_n)_{n\in {\BbbN}}$ $(\delta_n > 0)$
and constants $a_1, a_2 \in
{\BbbR}, b_1, b_2 > 0$ such that the following limits exist 
$$
\leqalignno{	
&\lim\limits_{n\to\infty}{1\over \delta_n} \left({n + \beta_n\over 2n
   + \alpha_n + \beta_n} - \varepsilon_n\right) ~=~ {a_1 \over 2}
& (2.4)\cr
   & &\cr
&\lim\limits_{n\to\infty} {1\over \delta_n} \left({n (n + \alpha_n) +
   (n + \beta_n) (n + \alpha_n + \beta_n)\over (2n + \alpha_n + \beta_n)^2}
   - \varepsilon_n\right) ~=~ {a_2\over 2} & (2.5)\cr
   & &\cr
&\lim\limits_{n\to\infty} {1\over \delta_n^2} {(n + \beta_n) (n +
   \alpha_n) n\over (2n + \alpha_n + \beta_n)^3} ~=~ {b_1\over 4}
& (2.6) \cr
   & &\cr
&\lim\limits_{n\to\infty} {1\over \delta_n^2} {(n + \beta_n) (n +
   \alpha_n) (n + \alpha_n + \beta_n) n\over (2n + \alpha_n + \beta_n)^4}
   ~=~ {b_2\over 4}& (2.7) \cr
}
$$
then
$$
\lim\limits_{n\to\infty} {1\over n} N^{(\alpha_n, \beta_n)} (\delta_n
   \xi + 2\varepsilon_n - 1) ~= \int_{a_2 - 2\sqrt{b_2}}^\xi
  ~ f^{(a_1, a_2, b_1, b_2)} (x) dx
\leqno (2.8)
$$
where the limiting density is given by
$$
f^{(a_1, a_2, b_1, b_2)} (x)~\colon = ~~~~~~~~~~~~~~~~~~~~~~~~~
~~~~~~~~~~~~~~~~~~~~~~~~~~~~~~~~~~~~
~~~~~~~~~~~~~~~~~~~~~~~~~~~~~~
$$
$$
   \cases{\displaystyle{
   {b_1\over 2\pi}{\sqrt{4b_2 - (x - a_2)^2}\over (b_2 - b_1)
             x^2 + (b_1 a_2 + b_1 a_1 - 2 b_2 a_1) x + b_2 a_1^2 -
		   a_1 a_2 b_1 + b_1^2}}  &if $|x - a_2| \le 2\sqrt{b_2}$\cr
		\noalign{\bigskip}
         ~ 0& ~~~~~~~~otherwise.\cr}
$$

\bigskip\noindent
{\bf Proof:}
Let $\mu_n$ denote the uniform distribution on the set
$$
\left\{{x - 2 \varepsilon_n + 1\over \delta_n}~ \bigg|~ P_n^{(\alpha_n,
   \beta_n)} (x) = 0\right\}
$$
then we obtain from (2.2) and Proposition 2.1 for the corresponding
Stieltjes transform
\newbox\DigitBox \setbox\DigitBox=\hbox{0}
\newdimen\DigitWidth \DigitWidth=\wd\DigitBox
\catcode`@=\active \def@{\kern\DigitWidth}
$$
\eqalign{
\Phi_{\mu_n} (z) ~=~ \int {d\mu_n (x)\over z - x} &~=~ \delta_n
   \int {d\mu^{(\alpha_n, \beta_n)} (x)\over \delta_n z +
   2\varepsilon_n - 1 - x}\cr
&~=~ {\lower5pt\hbox{@@@1@@@$|$}\over \raise5pt\hbox{$|z-\eta_1^{(n)}$}} -
   {\lower5pt\hbox{@@@$\lambda_1^{(n)}@|$} \over \raise5pt\hbox{$|z -
   \eta_2^{(n)}$}} - \ldots - {\lower5pt\hbox{$@\lambda_{n-1}^{(n)}@@|$}
   \over \raise5pt\hbox{$|z - \eta_n^{(n)}$}}\cr
} \leqno (2.9)
$$
where $p_0^{(n)} = 0,~ p_{-1}^{(n)} = 0$,
$$
\eqalign{
\lambda_j^{(n)} &~=~ {4\over \delta_n^2} (1 - p_{2j-2}^{(n)})
   p_{2j-1}^{(n)} (1 - p_{2j-1}^{(n)}) p_{2j}^{(n)} \qquad ~~(j = 1,
   \ldots, n-1),\cr
   & \cr
\eta_j^{(n)} &~=~ {2\over \delta_n} \displaystyle{
\left\{(1 - p_{2j-3}^{(n)}) p_{2j-2}^{(n)} +
   (1 - p_{2j-2}^{(n)}) p_{2j-1}^{(n)} - \varepsilon_n
   \right\}}~~~~~~ (j = 1,
   \ldots, n)\cr
}
$$
and the  $p_j^{(n)}$ are defined in (2.3).
From the limit assumptions (2.4) and (2.6)  of the theorem and
(2.3) it follows that  
$$
\lim_{n \to \infty} \lambda_1^{(n)}~=~\lim_{n\to \infty}
{4\over \delta_n^2} {(n-1)(n+\alpha_n)(n+\beta_n) \over
(2n+\alpha_n+\beta_n)^2(2n-1+\alpha_n+\beta_n)}~=~b_1
$$
and
$$
\lim_{n\to \infty} \eta_1^{(n)}~=~ {2\over \delta_n} \left(
{n+\beta_n \over 2n+\alpha_n+\beta_n} - \varepsilon_n \right)
~=~ a_1~~.
$$
Similary we have from (2.5) and (2.7) 
$$
\lim_{n\to\infty}
   \lambda_j^{(n)} = b_2~~,
 \lim_{n\to\infty}
   \eta_j^{(n)} = a_2 \ \ ~~~(j \ge 2).
$$
Because the quantities in the continued fraction expansion of the
Stieltjes transform are bijective continuous functions of the moments
of the given distribution we obtain that the moments of $\mu_n$
converge to the moments of a distribution $\mu$ with Stieltjes
transform
\newbox\DigitBox \setbox\DigitBox=\hbox{0}
\newdimen\DigitWidth \DigitWidth=\wd\DigitBox
\catcode`@=\active \def@{\kern\DigitWidth}
$$
\eqalign{
\Phi_\mu (z) &~=~ {\lower5pt\hbox{@@1@@$|$}\over \raise5pt\hbox{$|z-
   a_1$}} - {\lower5pt\hbox{@@$b_1@@|$}\over \raise5pt\hbox{$|z-
   a_2$}} - {\lower5pt\hbox{@@$b_2@@|$}\over \raise5pt\hbox{$|z-
   a_2$}} - {\lower5pt\hbox{@@$b_2@@|$}\over \raise5pt\hbox{$|z-
   a_2$}} - \ldots\cr
&~=~ {1\over 2} {(2b_2 - b_1) z + (b_1 a_2 - 2b_2 a_1) - b_1 \sqrt{(z -
   a_2)^2 - 4b_2}\over (b_2 - b_1) z^2 + z (b_1 a_2 + b_1 a_1 - 2b_2 a_1)
   + b_2 a_1^2 - a_1 a_2b_1  + b_1^2}\cr
}
$$
(here the last equality follows by straightforward algebra).
From Theorem 40 of Nevai ([12], p. 143)
we have that $\mu$ has a density in
the interval $[a_2 - 2 \sqrt{b_2}, a_2 + 2\sqrt{b_2}]$ and by the
inversion formula for Stieltjes transforms (see e.g.  [13])
this density is given by
$$
-{1\over \pi} \hbox{Im} (\Phi_\mu (x)) = f^{(a_1, a_2, b_1, b_2)} (x)\ \
\hbox{~~~~~~if ~}
   x\in [a_2 - 2\sqrt{b_2}, a_2 + 2 \sqrt{b_2}].
$$
Therefore the limit distribution $\mu$ is determined by its moments
and it follows from well known results of probability theory (see
e.g.~Feller [5], p.~263) that $\mu_n$ converges weakly with limit
$\mu$, that is
$$
\lim_{n\to\infty} {1\over n} N^{(\alpha_n, \beta_n)} (\delta_n \xi +
   2\varepsilon_n - 1) ~=~ \lim_{n\to\infty} \mu_n (\xi)
~ =~ \mu (\xi)~=
   \int_{a_2 - 2 \sqrt{b_2}}^\xi ~ f^{(a_1, a_2, b_1, b_2)} (x) dx
$$
\hfill \bull

\bigskip

In the following we will apply Theorem 2.2 to special sequences
$\alpha_n, \beta_n$ which might be of interest in applications.
All results can be proved by straightforward (but sometimes
tedious) calculations checking
the conditions of Theorem 2.2 and we only sketch the proof
of Corollary 2.4.
Corollary 2.3 has already been established in 
[7] as an application of strong asymptotic results for
the Jacobi polynomials, the other results are new.

\bigskip\noindent
{\bf Corollary 2.3.}~~
Let $\lim\limits_{n\to\infty} {\alpha_n\over n} = a$ and
$\lim\limits_{n\to\infty} {\beta_n\over n} = b \ (a, b \ge 0)$ then
$$
\lim_{n\to\infty} {1\over n} N^{(\alpha_n,\beta_n)}
 (\xi) ~=~ {2 + a + b\over 2\pi}
   \int_{r_1}^\xi {\sqrt{(r_2 - x) (x - r_1)}\over 1 - x^2} dx
   ~~~~~~~\quad r_1 \le
   \xi\le r_2
$$
where
$$
{r_{1,2}}~\colon =~ {b^2 - a^2 \pm 4 \sqrt{(a+1)(b+1)(a+b+1)}\over
(2+a+b)^2}~~.
$$

\bigskip\noindent
{\bf Corollary 2.4.}~~
Let $\lim\limits_{n\to\infty} {\alpha_n\over n} = \infty,
\lim\limits_{n\to\infty} {\beta_n\over n} = \infty$ and
$\lim\limits_{n\to\infty} {\alpha_n\over \beta_n} = c > 0$, then
$$
\lim_{n\to\infty} {1\over n} N^{(\alpha_n, \beta_n)}
   \left(\sqrt{{n\over \alpha_n}} \xi - {\alpha_n - \beta_n\over \alpha_n
   + \beta_n}\right) ~=~ {2\over \pi \sigma^2} \int_{-\sigma}^\xi
   \sqrt{\sigma^2 - x^2} dx \ \ ~~~~~|\xi| \le \sigma
$$
where $\sigma = {4c /(1 + c)^{3/2}}$.

\bigskip
\noindent
{\bf Proof:} In Theorem 2.2 we put $\delta_n= \sqrt{n/\alpha_n}$
and $\varepsilon _n=\beta_n/(\alpha_n+\beta_n)$. The asumptions
(2.4), (2.5), (2.6) and (2.7) are satifisfied with $a_1=a_2=0$,
$b_1=b_2=4c^2/(1+c)^3$ and the assertion follows from (2.8).
\hfill \bull

\bigskip\noindent
{\bf Corollary 2.5.}~~
Let $\lim\limits_{n\to\infty} {\alpha_n\over n} = \infty$ and
$\lim\limits_{n\to\infty} {\beta_n\over n} = b\ge 0$, then
$$
\lim_{n\to\infty} {1\over n} N^{(\alpha_n, \beta_n)} ({n\over \alpha_n} \xi
   - 1)~ =~ {1\over 4\pi} \int_{s_1}^\xi {\sqrt{(s_2-x) (x-s_1)}\over x}
 dx
~~~~~~~~~   \ s_1 \le \xi \le s_2
$$
where $s_{1,2} = 2 (2 + b) \pm 4 \sqrt{1 + b}$.

\bigskip\noindent
{\bf Corollary 2.6.}~~
Let $\lim\limits_{n\to\infty} {\alpha_n\over n} = \infty,
\lim\limits_{n\to\infty} {\beta_n\over n} = \infty$ and
$\lim\limits_{n\to\infty} {\alpha_n\over \beta_n} = \infty$, then
$$
\lim_{n\to\infty}{1\over n}
 N^{(\alpha_n, \beta_n)} \left({\sqrt{n\beta_n}\over
   \alpha_n} \xi - {\alpha_n + 2 \sqrt{n\beta_n} - \beta_n\over 2n +
   \alpha_n + \beta_n}\right) ~=~ {1\over 8\pi} \int_{-2}^\xi
   \sqrt{(6 - x) (x + 2)} dx
   $$
   for all $ -2 \le\xi\le 6$.

\bigskip\noindent
{\bf Example 2.7.}\ \ Letting $\alpha_n = n^4, \beta_n = n^3$ in
Corollary 2.6 we have that
$$
\lim_{n\to\infty}{1\over n}
 N^{(n^4, n^3)} \left({\xi\over n^2} - {n^3 - n^2 + 2n\over
   n^3 + n^2 + 2}\right) ~=~ {1\over 8\pi} \int_{-2}^\xi \sqrt{(6-x)
   (x+2)} \ dx
$$
for all $ -2 \le\xi\le 6$, or equivalently 
$$
\lim_{n\to\infty}{1\over n}
 N^{(n^4, n^3)} \left({4\over n^2} \xi + {2\over n^2} - {n^3
   - n^2 + 2n\over n^3 + n^2 + 2}\right) ~=~ {2\over \pi}
   \int_{-1}^\xi \sqrt{1 - x^2} dx.
\leqno(2.10)
$$ 
It should be noted that the location sequences $2\varepsilon_n - 1$ in
Theorem 2.2 can be changed by an addition of an amount $g_n$, and
maintain the same limit, only if $g_n/\delta_n \to 0$.
This means that the limit distribution of ${1\over n} N^{(n^4, n^3)}
({4\over n^2} \xi - 1)$ is not the same as (2.10).
We can however simplify this to
$$
\lim_{n\to\infty} {1\over n} N^{(n^4, n^3)} \left({4\xi\over n^2} - 1
   + {2n-2\over n^2}\right) ~=~ {2\over \pi} \int_{-1}^\xi \sqrt{1 -
   x^2} dx
$$
 \hfill\bull
 \bigskip

\noindent
{\bf Remark 2.8}\ \ The results of Corollary 2.3 -- 2.6 can be
motivated heuristically in the following way.
As an example consider the situation of Corollary 2.4 for $\alpha_n =
\beta_n$ (i.e.\ $c = 1$).
Because the polynomials $P_n^{(\alpha_n, \alpha_n)} (\sqrt{n/\alpha_n}
x)$ are orthogonal on the interval $[-\sqrt{\alpha_n/n},
\sqrt{\alpha_n/n}]$ with respect to the weight function
$$
\left(1 - {n\over \alpha_n} x^2\right)^{\alpha_n} \approx ~e^{-nx^2}
$$
we may expect that the uniform distribution on the set
$$
\left\{\sqrt{\alpha_n \over n} x~|~P_n^{(\alpha_n, \alpha_n)} (x) =
0\right\}
$$
has the same limit behaviour as the uniform distribution on the set
$$
\{{x \over \sqrt{n}}~|~H_n (x) = 0\}
$$
(here $H_n (x)$ denotes the $n$-th Hermite polynomial).
The weak limit of this distribution is known to be the measure with
density $\pi^{-1} \sqrt{2 - x^2}$ on the interval $[-\sqrt{2},
\sqrt{2}]$ (see e.g.\  [3]).

The case $\alpha_n \not= \beta_n$ is similar.
Thus the Jacobi weight function when properly scaled is approximately
the Hermite weight function, under the conditions of Corollary 2.4.
If $\alpha_n$ and $\beta_n$ go to infinity faster than $n$ we might
except the zeros of $P_n^{(\alpha_n \beta_n)} (x)$ to behave like the
Hermite polynomials when properly scaled.

 If $a$ and $b$ are both zero in Corollary 2.3 the limit density is the
classical arc-sin distribution.
So if $\alpha_n$ and $\beta_n$ both go to infinity slower than $n$ as
in Corollary 2.3, then the zero of $P_n^{(\alpha_n, \beta_n)} (x)$
behave like the arc-sin distribution.
 The limit distributions in Corollary 2.5 and 2.6 are somewhat similar.
The limit distribution in Corollary 2.5 is like the Laguerre case and
Corollary 2.6 is like the Hermite.
\hfill\bull

\bigskip
\bigskip\noindent
{\bf 3.\ \ The asymptotic distribution of the zeros of the generalized
Laguerre and Hermite polynomials.}
Throughout this section let
$$
\eqalign{
&N^{(\alpha_n)} (\xi)\colon ~=~ \#
\left\{~x~|~L_n^{(\alpha_n)} (x) = 0,~ x\le
   \xi\right\}\cr
&M^{(\gamma_n)} (\xi)\colon ~=~ \#
\left\{~x~|~H_n^{(\gamma_n)} (x) = 0,~ x\le
   \xi\right\}\cr
}
$$
denote the number of zeros of the generalized Laguerre and Hermite
polynomials less or equal than $\xi$, respectively.
The following theorem can be proved by similar arguments as in Section
2 using the corresponding characterization for the Laguerre
polynomials in  [3].
Its first part was also derived in  [6] as an application
of strong asymptotics for the generalized Laguerre polynomials and is
given here for the sake of completeness.

\bigskip\noindent
{\bf Theorem 3.1.}

\item{a)}Let $\lim\limits_{n\to\infty} {\alpha_n\over n} = a\ge 0$,
then
$$
\lim_{n\to\infty} {1\over n} N^{(\alpha_n)} (n\xi) ~=~ {1\over 2\pi}
   \int_{r_1}^\xi {\sqrt{(r_2-x) (x-r_1)}\over x} dx \qquad
   ~~~~r_1\le \xi\le r_2~,~
$$
where $r_{1,2} = 2 + a \pm 2\sqrt{1 + a}$
\medskip
\item{b)}Let $\lim\limits_{n\to\infty} {\alpha_n\over n} = \infty$,
then
$$
\lim_{n\to\infty} {1\over n} N^{(\alpha_n)} (\sqrt{n\alpha_n} \xi +
   \alpha_n) ~=~ {1\over 2\pi} \int_{-2}^\xi \sqrt{4 - x^2} dx
   \qquad ~~~~~~~|\xi| \le 2~~.
$$

\bigskip
\noindent
The corresponding results for the generalized Hermite polynomials can
easily be derived from Theorem 3.1 and the relations ($n \ge 0$)
$$
\eqalign{
H_{2n}^{(\gamma_n)} (x) &~=~ (-1)^n 2^{2n} n! L_n^{(\gamma_n - 1/2)}
   (x^2)\cr
H_{2n+1}^{(\gamma_n)} (x) &~=~ (-1)^n 2^{2n+1} n! x L_n^{(\gamma_n -
   1/2)} (x^2)\cr
}
$$
(see e.g.~Chihara [1]).

\bigskip\noindent
{\bf Theorem 3.2.}

\item{a)}Let $\lim\limits_{n\to\infty} {\gamma_n\over n} = c \ge 0$,
then
$$
\lim_{n \to \infty}
{1\over n} M^{(\gamma_n)} (\sqrt{n} \xi) = \int_{-\sigma}^\xi
   f (x) dx \qquad -\sigma \le \xi \le \sigma
$$
where $\sigma = \sqrt{1 + c + \sqrt{1 + 2c}}$,
$ \rho = \sqrt{1 + c -\sqrt{1 + 2c}}$ ~and
$$
f (x) ~=~ \cases{\displaystyle{
{1\over \pi} {\sqrt{(x^2 - \rho^2) (\sigma^2 - x^2)}
                \over |x|}} &
			    if $\rho < |x| < \sigma$\cr
			\noalign{\medskip}
               0& otherwise ~.\cr}
$$

\item{b)}Let $\lim\limits_{n\to\infty} {\gamma_n\over n} = \infty$,
then
$$
\lim_{n \to \infty}
{1\over n} M^{(\gamma_n)} \left(\sqrt{\sqrt{\gamma_n n} \xi^2 +
   \gamma_n}\right) = {1\over 2} + {1\over \pi} \int_0^\xi
   x \sqrt{2 - x^4} dx \qquad 0\le \xi\le \root 4 \of 2~.
$$

\bigskip\noindent
{\bf Remark 3.3.}~~
The first part of Theorem 3.2 can also be found in 
[6].
It is also worthwhile to mention that in contrary
to all other results in this paper the normalizing
sequence in Theorem 3.2b) is nonlinear. 
\hfill \bull

\bigskip
\bigskip\noindent
{\bf 4.\ \ Asymptotics for the largest and smallest zero.}
Using the Sturm Comparison Theory 
the following limit theorem for the largest and smallest zero of the
Jacobi polynomials $P_n^{(\alpha_n, \beta_n)} (x)$  was proved in [11]
(a different proof
can be found in a recent paper  [9] using general
results for bounds of the extreme zeros of orthogonal polynomials).

\bigskip\noindent
{\bf Theorem} (Moak, Saff, Varga) Let $r_n^{(\alpha_n,
\beta_n)}$ and $s_n^{(\alpha_n, \beta_n)}$ denote, respectively, the
smallest and largest zeros of the generalized Jacobi polynomials
$P_n^{(\alpha_n, \beta_n)} (x)$.
If
$$
\lim_{n\to\infty} {\alpha_n\over 2n + \alpha_n + \beta_n} = A \hbox{
  ~~ and ~~} \lim_{n\to\infty} {\beta_n\over 2n + \alpha_n + \beta_n} = B,
$$
then
$$
\lim_{n\to\infty} r_n^{(\alpha_n, \beta_n)} = r_{A, B}
   \hbox{ ~~and~~ } \lim_{n\to\infty} s_n^{(\alpha_n, \beta_n)} = s_{A, B}
$$
where
$$
r_{A, B},s_{A, B} = B^2 - A^2 \pm \sqrt{(A^2 + B^2 - 1)^2 - 4A^2 B^2}~~.
$$
Furthermore the zeros of the sequence $\{P_n^{(\alpha_n, \beta_n)}
(x)\}_{n=0}^\infty$ are dense in the interval $[r_{a, b}, s_{a, b}]$.

\bigskip

In the cases considered in Corollaries 2.4, 2.5 and 2.6 it can easily be
shown that $A+B = 1$ and the interval $[r_{A, B}, s_{A, B}]$
degenerates to the points ${1-c\over 1+c}, -1$ and $-1$, respectively.
The results of the previous sections
suggest that we can find similar statements as in the above theorem 
 when the sequences $\alpha_n$ and/or
$\beta_n$ converge to infinity with a larger order than $n$,
provided that the zeros of the orthogonal polynomials
are standardized appropriately. In fact we
have  the following second
order approximations for the zeros of the generalized Jacobi
polynomials $P_n^{(\alpha_n, \beta_n)} (x)$ when $\alpha_n/n \to
+\infty$ and/or $\beta_n/n \to +\infty$.

\bigskip\noindent
{\bf Theorem 4.1.}~~Let $r_n^{(\alpha_n,
\beta_n)}$ and $s_n^{(\alpha_n, \beta_n)}$ denote, respectively, the
smallest and largest zero of the generalized Jacobi polynomials
$P_n^{(\alpha_n, \beta_n)} (x)$, assume
that $\lim\limits_{n\to\infty} {\alpha_n\over n} = \infty,
\lim\limits_{n\to\infty} {\beta_n\over n} = \infty$ and
$\lim\limits_{n\to\infty} {\alpha_n\over \beta_n} = c > 0$, then
$$
\eqalign{
&\lim_{n\to\infty} \sqrt{{\alpha_n\over n}} \left\{r_n^{(\alpha_n,
   \beta_n)} + {\alpha_n - \beta_n\over \alpha_n + \beta_n}\right\}
  ~ =~ - {4c\over (1 + c)^{3/2}}\cr
  & \cr
&\lim_{n\to\infty} \sqrt{{\alpha_n\over n}} \left\{s_n^{(\alpha_n,
   \beta_n)} + {\alpha_n - \beta_n\over \alpha_n + \beta_n}\right\}
   ~=~ {4c\over (1 + c)^{3/2}}.\cr
}
$$
Furthermore, the zeros of
$$
P_n^{(\alpha_n, \beta_n)}\left(\sqrt{{n\over \alpha_n}}\left[x
   +{\beta_n - \alpha_n\over\beta_n + \alpha_n}\right]\right)
$$
become dense in the interval $[-4c (1 + c)^{-3/2}, 4c (1 + c)^{-3/2}]$.

\bigskip\noindent
{\bf Proof:} By Theorem 2 of  [9] we obtain for the
largest zero $s_n^{(\alpha_n, \beta_n)}$ of $P_n^{(\alpha_n, \beta_n)}
(x)$ the upper bound
$$
s_n^{(\alpha_n,\beta_n)} ~\le ~ \max \{ s_n(k)~|~ 0<k<n\}
\leqno (4.1)
$$
where 
$$
\eqalign{
&s_n(k) ~=~ 
{\beta_n^2 - \alpha_n^2\over (2k-2 + \beta_n + \alpha_n) (2k+2 + \beta_n +
   \alpha_n)} ~+~ {2\over 2k + \alpha_n + \beta_n}
\times \cr
& \cr
& \left[ \left({\beta_n^2 - \alpha_n^2 \over (2k -2
+\alpha_n +\beta_n)(2k+2 +\alpha_n +\beta_n)} \right)^2 +
{4k(k+\alpha_n)(k+\beta_n)(k+\alpha_n+\beta_n) \over
(2k-1+\alpha_n+\beta_n)(2k+1+ \alpha_n+\beta_n)}\right]^{1\over 2}\cr
}
$$
(note that there is a typographical error in the paper 
 [9], p. 138). If $\beta_n \ge \alpha_n$ 
it follows that  
$$
s_n(k) ~\le ~ s_n^{(1)}~:= ~
{\beta_n - \alpha_n\over \beta_n +
   \alpha_n} ~+~ {2\over 2 + \alpha_n + \beta_n} \sqrt{g (\alpha_n,
   \beta_n) + h (\alpha_n, \beta_n)}
\leqno(4.2)
$$
where
$$
g (\alpha_n, \beta_n) :=  \left({\beta_n - \alpha_n\over \beta_n +
   \alpha_n}\right)^2
$$
and
$$
h (\alpha_n, \beta_n) := {4 n (n + \alpha_n) (n + \beta_n) (n +
   \alpha_n + \beta_n)\over (\alpha_n + \beta_n + 1) (\alpha_n + \beta_n
   + 3)}.
$$
On the other hand, if $\beta_n \le \alpha_n$, we have
$$
s_n(k) ~\le ~ s_n^{(2)}~:= ~
{\beta_n^2 - \alpha_n^2\over (\beta_n +
   \alpha_n+2n)^2} ~+~ {2\over 2 + \alpha_n + \beta_n} \sqrt{g (\alpha_n,
   \beta_n) + h (\alpha_n, \beta_n)}
\leqno(4.3)
$$
and consequently from (4.1), (4.2) and (4.3)
$$
s_n^{(\alpha_n,\beta_n)} ~\le ~ \max \{ s_n^{(1)}, s_n^{(2)} \}~.
\leqno (4.4)
$$
Now 
straightforward algebra now yields that
$$
\leqalignno{
&\lim_{n\to\infty} {\alpha_n\beta_n
\over (2 + \alpha_n + \beta_n)^2} \cdot g
   (\alpha_n, \beta_n) ~=~ {c (c-1)^2\over (c+1)^4}~~~, &(4.5)\cr
&\lim_{n\to\infty} {\alpha_n\over n} {1\over (2 + \alpha_n +
   \beta_n)^2}  \cdot
   h (\alpha_n, \beta_n) ~=~ {4c^2\over (1+c)^3}~~~, &(4.6)\cr
}
$$
and  by a combination  of (4.2), (4.3), (4.5) and (4.6) it follows that  
$$
\mathop{\lim \sup}\limits_{n \to \infty} \sqrt{\alpha_n \over n}
\left\{ s_n^{(j)} +{\alpha_n -\beta_n \over 
\alpha_n +\beta_n} \right) ~\le ~
{4c \over (1+c)^{3/2}}~~~~~~~~~~~~~j=1,2~.
$$
This implies (observing (4.4))
$$
\mathop{\lim\ \sup}\limits_{n\to\infty} \sqrt{{\alpha_n\over n}}
   \left\{s_n^{(\alpha_n, \beta_n)} + {\alpha_n - \beta_n\over \alpha_n +
   \beta_n}\right\} ~\le ~{4c\over (1 + c)^{3/2}}.
\leqno(4.7)
$$
\noindent
It is clear from Corollary 2.4 that equality must hold in (4.7); which
proves the assertion regarding the largest zero.
The corresponding statement for the smallest zero $r_n^{(\alpha_n,
\beta_n)}$ follows from the well known relation $P_n^{(\alpha, \beta)}
(x) = P_n^{(\beta, \alpha)} (-x)$ and the first part.
Finally the statement regarding the denseness of the zeros also
follows from Corollary 2.4. \hfill \bull

\bigskip
\noindent
In the following theorems we consider the remaining cases of Corollary
2.5 and 2.6 for the Jacobi polynomials and the corresponding results
for the generalized Laguerre and Hermite polynomials.
The proofs are similar to that of Theorem 4.1 and therefore omitted.

\bigskip\noindent
{\bf Theorem 4.2.}~~Let $r_n^{(\alpha_n,
\beta_n)}$ and $s_n^{(\alpha_n, \beta_n)}$ denote, respectively, the
smallest and largest zero of the generalized Jacobi polynomials
$P_n^{(\alpha_n, \beta_n)} (x)$ and
assume that $\lim\limits_{n\to\infty} {\alpha_n\over n} = \infty$ and
$\lim\limits_{n\to\infty} {\beta_n\over n} = b \ge 0$, then
$$
\eqalign{
\lim_{n\to\infty} {\alpha_n\over n} \{s_n^{(\alpha_n, \beta_n)} + 1\}
   &~=~ 2 (2 + b) + 4\sqrt{1 + b}\cr
   &\cr
\lim_{n\to\infty} {\alpha_n\over n} \{r_n^{(\alpha_n, \beta_n)} + 1\}
   &~= ~2 (2 + b) - 4\sqrt{1 + b}~~.\cr
}
$$
Furthermore, the zeros of
$$
P_n^{(\alpha_n, \beta_n)} \left({n\over \alpha_n} [x - 1]\right)
$$
become dense in the interval $[2 (2+b) - 4 \sqrt{1+b}, 2 (2+b) + 4
\sqrt{1+b}]$.

\bigskip\noindent
{\bf Theorem 4.3.}~~ Let $r_n^{(\alpha_n,
\beta_n)}$ and $s_n^{(\alpha_n, \beta_n)}$ denote, respectively, the
smallest and largest zero of the generalized Jacobi polynomials
$P_n^{(\alpha_n, \beta_n)} (x)$, assume that $\lim\limits_{n\to\infty}
{\alpha_n\over n} = \infty$, $\lim\limits_{n\to\infty} {\beta_n\over n}
= \infty $, $\lim\limits_{n\to\infty} {\alpha_n\over \beta_n} = \infty$
and let $\varepsilon_n = (\alpha_n + 2 \sqrt{n\beta_n} - \beta_n) /
(2n + \alpha_n + \beta_n)$, then
$$
\eqalign{
&\lim_{n\to\infty} {\alpha_n\over \sqrt{n\beta_n}} \{s_n^{(\alpha_n,
   \beta_n)} + \varepsilon_n\} ~=~ 6\cr
& \cr
&\lim_{n\to\infty} {\alpha_n\over \sqrt{n\beta_n}} \{r_n^{(\alpha_n,
   \beta_n)} + \varepsilon_n\}~ =~ -2.\cr
}
$$
Furthermore, the zeros of the sequence
$$
P_n^{(\alpha_n, \beta_n)} \left({\sqrt{n\beta_n}\over \alpha_n}
   [x - \varepsilon_n]\right)
$$
become dense in the interval $[-2, 6]$.

\bigskip\noindent
{\bf Theorem 4.4.}~~ Let $r_n^{(\alpha_n)}$ and $s_n^{(\alpha_n)}$,
respectively, denote the smallest and largest zero of the generalized
Laguerre polynomials $L_n^{(\alpha_n)} (x)$.

\item{a)}Assume that $\lim\limits_{n\to\infty} {\alpha_n\over n} = a
\ge 0$, then
$$
\eqalign{
\lim_{n\to\infty} {s_n^{(\alpha_n)}\over n} &~=~ 2 + a + 2 \sqrt{1+a}\cr
\lim_{n\to\infty} {r_n^{(\alpha_n)}\over n} &~=~ 2 + a - 2 \sqrt{1+a}~~.\cr
}
$$
Furthermore, the zeros of $L_n^{(\alpha_n)} (nx)$ are dense in the
interval $[2 + a - 2 \sqrt{1+a}, 2 + a + 2 \sqrt{1+a}]$.

\item{b)}Assume that $\lim\limits_{n\to\infty} {\alpha_n\over n} =
\infty$, then
$$
\eqalign{
&\lim_{n\to\infty} {s_n^{(\alpha_n, \beta_n)} - \alpha_n\over \sqrt{n
   \alpha_n}}~ =~ 2\cr
&\lim_{n\to\infty} {r_n^{(\alpha_n, \beta_n)} - \alpha_n\over \sqrt{n
   \alpha_n}}~ =~ -2~~.\cr
}
$$
Furthermore, the zeros of $L_n^{(\alpha_n)} (\sqrt{n \alpha_n} x +
\alpha_n)$ become dense in the interval $[-2, 2]$.

\bigskip\noindent
{\bf Theorem 4.5.}~~ Let $\rho_n^+ ~(\rho_n^-)$ and $\sigma_n^+
~(\sigma_n^-)$ denote, respectively, the smallest and largest
positive (negative) zero of the generalized Hermite polynomials
$H_n^{(\gamma_n)} (x)$.

\item{a)}Assume that $\lim\limits_{n\to\infty} {\gamma_n\over n} = c
\ge 0$, then
$$
\eqalign{
\lim_{n\to\infty} {\sigma_n^+ \over \sqrt{n}} &~=~ -\lim_{n\to\infty}
   {\rho_n^-\over \sqrt{n}} ~=~ \sqrt{1+c + \sqrt{1+2c}} ~=~ \sigma\cr
\lim_{n\to\infty} {\rho_n^+ \over \sqrt{n}} &~=~ -\lim_{n\to\infty}
   {\sigma_n^-\over \sqrt{n}} ~=~ \sqrt{1+c - \sqrt{1+2c}} ~=~ \rho ~~.\cr
}
$$
Furthermore, the zeros of $H_n^{(\gamma_n)} ({\sqrt{n}x})$ become
dense in $[-\sigma, -\rho] \cup [\rho,
\sigma]$.

\item{b)}Assume that $\lim\limits_{n\to\infty} {\gamma_n\over n} =
\infty$, then
$$
\eqalign{
&\lim_{n\to\infty} {(\sigma_n^+)^2 - \gamma_n\over \sqrt{n
   \gamma_n}} ~=~ \lim_{n\to\infty} {(\rho_n^-)^2 - \gamma_n\over
   \sqrt{n\gamma_n}} ~=~ \sqrt{2}\cr
&\lim_{n\to\infty} {(\sigma_n^-)^2 - \gamma_n\over \sqrt{n
   \gamma_n}} ~=~ \lim_{n\to\infty} {(\rho_n^-)^2 - \gamma_n\over
   \sqrt{n\gamma_n}} ~=~ 0~~.\cr
}
$$
\hfuzz=16pt
Furthermore, the zeros of $H_n^{(\gamma_n)}
 (\sqrt{\sqrt{n \gamma_n} x^2 + \gamma_n})$
become dense in the interval $[-\sqrt{2}, \sqrt{2}]$.

\bigskip
\medskip\noindent
{\bf Acknowledgements.}
The authors would like to thank Professor Dr.~W.~Gawronski and Dr.~C.
Bosbach for their interest and helpful comments during the preparation
of the paper.
Professor Gawronski has obtained similar asymptotic distribution
results using the characterizing differential equation for the
classical orthogonal polynomials.
A hint of Dr.~Bosbach led our attention to the paper of Ismail and Li
[9] which turned out to be useful for the proofs of the
results in Secton 4. We are also indebted to two unknown referees
for their helpful comments.

\medskip
\bigskip
\noindent
{\bf References:}
\advance \leftskip by 0.2in \parindent =-0.2in

{\baselineskip=11pt
\parskip=7pt

1. Chihara, T.S. (1978).
{\it Introduction to Orthogonal Polynomials}, Gordon \& Breach, New
York.

2. Chen, L.C. and Ismail, M.E.H. (1991).
On asymptotics of Jacobi polynomials, {\it SIAM J.\ Math.\ Anal}.,
{\bf 22}, 1442--1449.

3. Dette, H. and Studden, W.J. (1992).
On a new characterization of the classical
orthogonal polynomials, {\it J. Approx. Theory}.,
{\bf 71}, 3--17.

4. Elbert, A. and Laforgia, A. (1987).
New properties of the zeros of a Jacobi polynomial in relation to
their centroid, {\it SIAM J.\ Math.\ Anal}., {\bf 18}, 1563--1572.

5. Feller, W. (1966).
{\it An Introduction to Probability Theory and Its Applications,
Vol.~II}, Wiley, New York.

6. Gawronski, W. (1993).
Strong asymptotics and the asymptotic zero distributions  of Laguerre
polynomials $L_n^{(an + \alpha)} (x)$ and Hermite polynomials
$H_n^{(an + \alpha)} (x)$, to appear: Analysis.

7. Gawronski, W. and Shawyer, B. (1991).
Strong asymptotics and the limit distribution of the zeros of the
Jacobi polynomials $P_n^{(an + \alpha, bn + \beta)}$, {\it Progress
in Approximation Theory}, 379--404, Academic Press, New York.

8. Gonchar, A.A. and Rakhmanov, E.A. (1986).
Equilibrium measure and the distribution of zeros of extremal
polynomials, {\it Math.\ USSR Sbornik}, {\bf 53}, 119--130.

9. Ismail, M.E.H. and Li, X. (1992).
Bound on the extreme zeros of orthogonal polynomials, {\it Proceedings
of the American Mathematical Society}, {\bf 115}, 131--140.

\advance \leftskip by 0.21in \parindent -0.2in
10. Mhaskar, H.N. and Saff, E.B. (1984).
Weighted polynomials on finite and infinite intervals: a unified
approach, {\it Bull.\ Amer.\ Math.\ Soc}., {\bf 11}, 351--354.

11. Moak, D., Saff, E.B. and Varga, R. (1979).
On the zeros of Jacobi polynomials $P_n^{(\alpha_n, \beta_n)} (x)$,
{\it Trans.\ Amer.\ Math.\ Soc}., {\bf 249}, 159--162.

12. Nevai, P. (1979).
{\it Orthogonal polynomials}, Memoirs Amer.\ Math.\ Soc., Vol.~213.

13. Perron, O. (1954).
{\it Die Lehre von den Kettenbr\"uchen (Band I and II)}, B.G. Teubner,
Stuttgart.

14. Skibinsky, M. (1986).
Principal representations and canonical moment sequences for
distributions on a interval, {\it J.\ Math.\ Anal.\ Appl}., {\bf 120},
95--120.

15. Szeg\"o, G. (1975).
Orthogonal Polynomials, {\it American Mathematical Society}
Colloquium Publications, Vol. {\bf 23}, Amer. Math. Soc., Providence,
RI.

16. Wall, H.S. (1948).
{\it Analytic theory of continued fractions}, Van Nostrand, New York.

}
\bye